\input amssym.tex
\def\1{\'{\i}} 

\def\Z{\hbox{$\Bbb Z$}}  
\def\N{\hbox{$\Bbb N$}}  
\def\C{\hbox{$\Bbb C$}} 
\def\qed{\vrule height6pt width4pt depth0pt}

\centerline{\bf Dirichlet series for complex powers of the Riemann zeta function}
\medskip
\rightline{Winston Alarc\'on Athens} 
\rightline{Universidad de Costa Rica}
\medskip
{\narrower\narrower\medskip\noindent{\bf Abstract.}
To obtain the Dirichlet series for complex powers of the Riemann zeta function, we define and study the basic properties of a sequence
 of polynomials
  that, used as coefficients of the respective terms of the Dirichlet series of the Riemann zeta function
  in the half plane $ x > 1 $, produces the required exponential function. Unlike the method described in ([4], p.~278), which requires more advanced knowledge of the relationships between Dirichlet series and multiplicative arithmetic functions, our approach only needs mathematical induction on the total number of prime divisors of $ n $, the Dirichlet product and the use of an analytic property characteristic of the exponential function in the complex plane.\medskip}
\bigskip
\parindent=0pt
{\bf I.\ Introduction.} 
\medskip
The problem to find the Dirichlet series for the powers of the function $ \zeta(s) = \sum_{n = 1}^\infty n^{- s} $ 
has been studied at least since the mid-nineteenth century. The first reference the author found is a work by V.~Bouniakowsky published in St.~Petersburg in {\oldstyle1862} ([1], p.~287), in which he considers the series 
$\psi(x) = \sum_{n = 1}^\infty {1 / n^x}$ and $ \big(\psi(x) \big)^m = \sum_{n = 1}^\infty {z_{n, m} / n^x} $ and then shows $ z_{n, 2} $ is the number $ N_0(n) $ of divisors of $ n $  and  $ z_{n, m} = N_{m-2}(n) = \sum_{d | n} N_{m-3}(d) $ for $ m = 3,4,5, \ldots$
\medskip 
Several authors ([2], p.~334; [5], p.~44-45) have obtained this result and have
identified the coefficient of the $n$-th term of the Dirichlet series for the $m$-th power of $\zeta(s)$ 
as the number of ways to express $n$ as the product of $m$ positive integers, taking into account the order of the factors (in Hardy \& Wright's description, denoted as $d_m(n)$); or, equivalently, as the number of solutions of the indeterminate equation $ x_1 \, x_2 \, \cdots \, x_m = n $, where the $ x_i \in \Z^+$ (in Vinogr\'adov's description, denoted as $\tau_m(n)$). They demonstrate the aforementioned arithmetic functions are multiplicative, so they limit themselves to proof the formula 
$$ 
d_m(n) = \tau_m(n) = {m + \nu-1 \choose m-1}
$$ 
for the case $ n = p^\nu $, with $ p $ any prime number and $m, \nu$ any positive integers.
\medskip 
The formula for the Dirichlet series of the {\it complex powers $z$\/} of the Riemann zeta function $\zeta(s)$ can be easily anticipated from the previous results, noting that 
$$ 
{m + \nu-1  \choose m-1} = {m + \nu-1 \choose \nu}
$$ 
and that the binomial coefficient ${w\choose \nu}$ can be generalized for any $w\in\C$ and $\nu\in\N$ by 
$$
{w\choose\nu}:={1\over\nu!} \prod_{0\leq j<\nu} (w-j) \eqno(1)
$$
In G\'erald Tenenbaum's {\it Introduction to Analytic and Probabilistic Number Theory\/} ([4], p.~278), the general formula for the complex powers of the Riemann zeta function is utilized in the study of the Delange -- Selberg method exposed in chapter II.5. The generalization for complex exponents $z$ is based on the properties of the Dirichlet series associated with the multiplicative arithmetic functions exposed in the chapter II.1.
\medskip
In what follows we obtain the general formula, considering the sequence of polynomials in the indeterminate $z$ defined by:
$$
\alpha_1(z) := 1, \quad 
\alpha_n(z) := \prod_{j\in{\Bbb Z}^+} {z+\nu_j -1 \choose \nu_j } \qquad n = 2,3,4,\ldots \eqno(2)
$$
where the $ \nu_j $ are the non-negative integer exponents in the canonical prime decomposition of $n$ of the form 
$$ 
n=\prod_{j \in {\Bbb Z}^+} p_j^{\nu_j} 
$$ 
$p_j$ denote the $j$-th prime number. 
\medskip
{\bf Examples.} For all complex numbers $ z $ and three distincts prime numbers $p$, $q$, $r$, we have:
$$
\alpha_p(z) = z,
$$
$$
\alpha_{p^2}(z)= {z(z+1)\over2},\qquad 
\alpha_{pq}(z)= z^2,$$
$$\alpha_{p^3}(z)= {z(z+1)(z+2)\over6},\qquad
\alpha_{p^2 q}(z)= {z^2(z+1)\over2},\qquad 
\alpha_{pqr}(z)=z^3.
$$
\bigskip
{\bf 2. Basic properties}.
\medskip
{\bf Theorem 1.} 
\smallskip
{\bf 1.1}\quad$\alpha_n(0) = 0$ for all integers $n>1$.
\smallskip
{\bf 1.2}\quad$\alpha_n(1)=1$ for all integers $n\geq 1$.
\smallskip
{\bf 1.3}\quad For each complex number $ z $, the function 
$n\mapsto \alpha_n(z)$ from $\Z^+$ to \C\ 
is a {\it multiplicative arithmetic function}; that is, provided $m$ and $n$ are coprime positive integers, the following equation is satisfied:
$$
\alpha_m(z) \, \alpha_n(z) = \alpha_{m n}(z).
$$
{\bf 1.4}\quad $$\sum_{j=0}^{k+1}\alpha_{p^j}(s)\,\alpha_{p^{k-j+1}}(t) = {s+t+k \over k+1} \sum_{j=0}^{k}\alpha_{p^j}(s)\,\alpha_{p^{k-j}}(t)\eqno(3)$$ for any pair $(s,t)$ of complex numbers, for all prime numbers $p$ and any integer $k\geq0$.
\smallskip
{\bf Proof.} 
The properties {\bf 1.1},  {\bf 1.2} and {\bf 1.3} are an immediate consequence of the definition equalities (1) and (2). For proof of (3), we can write:
$$\eqalign{
 {s+t+k \over k+1} \sum_{j=0}^{k} \alpha_{p^j}(s)\,\alpha_{p^{k-j}}(t) &=
  {1 \over k+1} \sum_{j=0}^{k} (s+t+k)\, \alpha_{p^j}(s)\,\alpha_{p^{k-j}}(t) 
  \cr 
&=  {1 \over k+1} \sum_{j=0}^{k} \big((s+j)+(t+k-j)\big)\, \alpha_{p^j}(s)\,\alpha_{p^{k-j}}(t) 
  \cr 
  &= {1 \over k+1} \sum_{j=0}^{k} (s+j)\, \alpha_{p^j}(s)\,\alpha_{p^{k-j}}(t) + (t+k-j)\,\alpha_{p^j}(s)\,\alpha_{p^{k-j}}(t) 
  \cr
&= {1 \over k+1} \Bigg( \sum_{j=0}^{k} (j+1)\, \alpha_{p^{j+1}}(s) \, \alpha_{p^{k-j}}(t) 
+ \sum_{j=0}^{k} (k-j+1)\,\alpha_{p^j}(s)\,\alpha_{p^{k-j+1}}(t) \Bigg)
}$$
since 
$$
(s+j)\,\alpha_{p^j}(s) = (s+j) {s+j-1\choose j} = (j+1)\,{s+j \choose j+1}=(j+1)\,\alpha_{p^{j+1}}(s)
$$ 
and 
$$
(t+k-j)\,\alpha_{p^{k-j}}(t) = (t+k-j)\,{t+k-j-1 \choose k-j} = (k-j+1)\,{t+k-j \choose k-j+1} = (k-j+1)\,\alpha_{p^{k-j+1}}(t).
$$
\smallskip 
Changing the summation index $ j $ to $ j-1 $ in the first summation on the right side of the last equation: \phantom{$\|$}
$$\eqalign{{s+t+k \over k+1} \sum_{j=0}^{k} \alpha_{p^j}(s)\,\alpha_{p^{k-j}}(t) &= {1\over k+1} \Bigg(\sum_{j=1}^{k+1} j\,\alpha_{p^j}(s)\,\alpha_{p^{k-j+1}}(t)
+ \sum_{j=0}^{k} (k-j+1)\,\alpha_{p^j}(s)\,\alpha_{p^{k-j+1}}(t) \Bigg) \cr
&= \alpha_{p^0}(s)\,\alpha_{p^{k+1}}(t)  
+ {1\over k+1} \sum_{j=1}^k (k+1)\,\alpha_{p^j}(s)\,\alpha_{p^{k-j+1}}(t)
+ \alpha_{p^{k+1}}(s)\,\alpha_{p^0}(t)
\cr
&= \sum_{j=0}^{k+1} \alpha_{p^j}(s)\,\alpha_{p^{k-j+1}}(t). \ \qed 
}$$
{\bf Theorem 2.} For any pair $ s $, $ t $ of complex numbers, for every integer $ n \geq 1 $, the following equality holds:
$$
\sum_{j|n} \alpha_j(s) \, \alpha_{n/j}(t) = \alpha_n(s + t)\eqno(4)
$$
where $ j | n $ under the summation symbol indicates the summation index $ j $ run over the set of all positive divisors of $ n $.
\smallskip
{\bf Proof.} We will first prove that (4) is true when $ n = p^k $, $p$ any prime number and $k$ any integer ${}\geq0$. That is, we will prove:
$$
\sum_{j=0}^{k}\alpha_{p^j}(s)\,\alpha_{p^{k-j}}(t) = \alpha_{p^k}(s+t). \eqno(5)
$$
We will prove (5) by induction on $k = 0,1,2,\ldots$\ 
\smallskip
For $ k = 0 $, by Theorem 1.2, (5) is trivially true.\ 
\smallskip
Assuming (5) is true for an arbitrary but fixed non negative integer $k$, using (3), we can write:
$$
\sum_{j=0}^{k+1}\alpha_{p^j}(s)\,\alpha_{p^{k-j+1}}(t) = {s+t+k \over k+1} \sum_{j=0}^{k}\alpha_{p^j}(s)\,\alpha_{p^{k-j}}(t) = {s+t+k \over k+1}\,\alpha_{p^k}(s+t) = \alpha_{p^{k+1}}(s+t) 
$$
which completes the proof of (5).
\medskip
{\bf Proof of (4)}. By induction on $\Omega(n)$, the total number of prime divisors of $n$, each counted according to its multiplicity, with $\Omega(1)=0$.
\smallskip
If $\Omega(n) = 0$, then $ n = 1 $ and, by Theorem 1.2, (4) is immediate.
\smallskip
Let $g \in \{0,1,2,\ldots\}$ be arbitrary but fixed and suppose that (4) is true for any $n \in \Z^+$ such that $\Omega(n) \leq g$.
\smallskip
Let $ n' $ be any positive integer such that $\Omega(n') = g + 1$. 
To avoid the case (5) that we have already proved, we can assume 
$ n'$ has at least two distinct prime divisors $ p $ and $ q $. Let $ k \in \Z^+ $ be such that $ p^k $ divides $ n'$ but $ p^{k + 1} $ does not divide $ n' $.
Let $ N = n' / p^k $. Since $ k \geq1 $, we have
$ \Omega(N) \leq g $,
so, by the induction hypothesis, (4) is true for this $N$.
\smallskip
We can partition the set of divisors of $ n'= p^{k} N $ into the following mutually disjoint subsets: the set $ A $ of divisors of $ N $ and the sets $ B_i = p^i A $ that result from multiplying each element of $ A $ by $ p^i $, where $ i $ runs through the set $ \{1, \ldots, k \} $. In this way, we can write:
$$
\sum_{j|n'} \alpha_j(s) \, \alpha_{n'/j}(t) 
= \sum_{j\in A} \alpha_j(s) \, \alpha_{p^k N/j}(t) + \sum_{i=1}^{k} \sum_{j\in B_i} \alpha_j(s) \, \alpha_{p^k N/j}(t) .
$$
Next, changing the summation index $ j $ to $ p^i j $ and $ B_i $ to $ A $ in the rightmost summation, we obtain:
$$\sum_{j|n'} \alpha_j(s) \, \alpha_{n'/j}(t) 
= \sum_{j\in A} \alpha_j(s) \, \alpha_{p^k N/j}(t) + \sum_{i=1}^{k} \sum_{j\in A} \alpha_{p^i j}(s) \, \alpha_{p^{k-i} N/j}(t).
$$
Now, for each $ j \in A $, $ p^k $ and $ N / j $ are coprime numbers. Also, for each $ i \in \{1,2, \ldots, k \} $ and each $ j \in A $, $ p^{i} $ and $ j $, as well as $ p^{k - i} $ and $ N / j $ are also pairs of coprime numbers. By Theorem 1.3, the hypothesis of induction and equation (5), we can write:
$$\eqalign{
\sum_{j|n'} \alpha_j(s) \, \alpha_{n'/j}(t)
&= \sum_{j\in A} \alpha_j(s) \, \alpha_{N/j}(t) \, \alpha_{p^k}(t) + \sum_{i=1}^{k} \sum_{j\in A} \alpha_{p^i}(s) \, \alpha_{j}(s) \, \alpha_{p^{k-i}}(t) \, \alpha_{N/j}(t) \cr
&= \alpha_{p^k}(t) \sum_{j \in A} \alpha_j(s) \, \alpha_{N/j}(t) + \sum_{i=1}^{k} \alpha_{p^i}(s)\,\alpha_{p^{k-i}}(t) \sum_{j \in A} \alpha_j(s) \, \alpha_{N/j}(t) \cr
&=  \Big(\alpha_{p^k}(t) + \sum_{i=1}^{k}\alpha_{p^i}(s)\,\alpha_{p^{k-i}}(t)\Big) \sum_{j | N} \alpha_j(s) \, \alpha_{N/j}(t) \cr
&=  \Big(\sum_{i=0}^{k}\alpha_{p^i}(s)\,\alpha_{p^{k-i}}(t)\Big)\, \alpha_N(s+t) \cr
&= \alpha_{p^k} (s+t) \, \alpha_N(s+t) \cr
&= \alpha_{n'}(s+t).\ \qed
}$$
\bigskip
{\bf 3. The general formula. }
\medskip
{\bf Theorem 3.} For each complex number $s$ whose real part is greater than unity, the function $\varphi_s: \C \to \C$ defined by
$$
\varphi_s(z) = \sum_{n=1}^\infty {\alpha_n(z) \over n^s} \eqno(6)
$$
is an exponential function in base $\zeta (s)$, that is:
$$ 
\sum_{n=1}^\infty {\alpha_n(z) \over n^s} = \exp\big(z \ln\zeta(s)\big) = \zeta(s)^z\eqno(7)
$$
{\bf Proof.} First note that, by Theorem 1.2, (7) holds for $z=1$. Now, 
since $\varphi_s: \C \to \C$ is an analytic function and $\varphi_s(1) = \zeta(s) \neq 0$ for all $s \in \C$ such that $\Re(s) > 1$, it is enough to show ([3], p.~110) that $\varphi_s$ satisfies the functional equation
$$
\varphi_s(z_1+z_2) = \varphi_s(z_1)\,\varphi_s(z_2). \eqno(8)
$$
The equality (8) is an immediate consequence of Theorem {\bf 2}, since we can apply the Dirichlet product to the series in the right side of equality (6): 
$$\eqalign{
\bigg(\sum_{n=1}^\infty \alpha_n(z_1) / n^s\bigg) \, \bigg(\sum_{n=1}^\infty \alpha_n(z_2) / n^s\bigg) 
&= \sum_{n=1}^\infty \bigg(\sum_{j|n} \alpha_j(z_1)\,\alpha_{n/j}(z_2)\bigg) / n^s \cr
&= \sum_{n=1}^\infty \alpha_n(z_1+z_2) / n^s \cr
&= \varphi_s(z_1+z_2).\ \qed
}$$
\bigskip
{\bf References.}
\medskip
[1]\  Dickson, Leonard E.\ {\it History of the Theory of Numbers}, Vol.~I, 3th ed. (1923).   
Reprint (1971): Chelsea.
\smallskip
[2]\  Hardy G.~H., Wright E.~M.\ {\it An Introduction to the Theory of Numbers}.  (1938). Heath-Brown, D.~R.; Silverman, J.~H.\ (eds.), 6th ed.\  (2008): Oxford University Press.
\smallskip
[3]\ Markushevich, A.~I. {\it Teor\1a de las funciones anal\1ticas}, Vol.~I.\ (Spanish trans.\ from Russian) (1987).\ Moscow, 2nd.\ repr.:  (1970). MIR.
\smallskip
[4]\ Tenenbaum, G.\ {\it Introduction to Analytic and Probabilistic Number Theory}. (English trans. from French). Third ed. (2008). Graduate Studies in Mathematics, Vol.\ 163. A.\ M.\ M.   
\smallskip
[5]\ Vinogr\'adov, I.\ M.\ {\it Fundamentos de la Teor\1a de los n\'umeros}.\ (Spanish trans.\ from Russian) (1987).\ Moscow, 2nd.\ repr. (1971): MIR.

\bye